\documentclass{article}
\usepackage[utf8]{inputenc}

\usepackage{amsmath}
\usepackage{amsfonts}
\usepackage{caption}
\usepackage{graphicx}
\usepackage{cite}
\usepackage{url}
\usepackage{booktabs}
\usepackage{bm}
\newcommand{\tikzcmark}{
\tikz[scale=0.23] {
    \draw[line width=0.7,line cap=round] (0.25,0) to [bend left=10] (1,1);
    \draw[line width=0.8,line cap=round] (0,0.35) to [bend right=1] (0.23,0);
}}

\usepackage{geometry}
\usepackage{tabularx,colortbl, makecell, caption, multirow}

\usepackage{todonotes}

\usepackage{xcolor}

\title{A computational study of off-the-shelf MINLP solvers on a benchmark set of
congested capacitated facility location problems}
\author{Pasquale Avella$^1$, Alice Calamita*$^2$, Laura Palagi$^2$\\\\
$^1$\small{DING, Università del Sannio, Benevento, 82100, Italy}\\\small{\texttt{avella@unisannio.it}}\\
$^2$\small{DIAG, Sapienza University of Rome, Rome, 00185, Italy}\\
\small{\texttt{\{alice.calamita,laura.palagi\}@uniroma1.it}}\\\small{* Corresponding author}}
\date{}

\begin{document}

\maketitle

\abstract 

This paper analyzes the performance of five well-known off-the-shelf optimization solvers on a set of congested capacitated facility location problems formulated as mixed-integer conic programs (MICPs).
We aim to compare the computational efficiency of the solvers and examine the solution strategies they adopt when solving instances with different sizes and complexity.
The solvers we compare are Gurobi, Cplex, Mosek, Xpress, and Scip. We run extensive numerical tests on a testbed of 30 instances from the literature.
Our results show that Mosek and Gurobi are the most competitive solvers, as they achieve better time and gap performance, solving most instances within the time limit. Mosek outperforms Gurobi in large-size problems and provides more accurate solutions in terms of feasibility. Xpress solves to optimality about half of the instances tested within the time limit, and in this half, it achieves performance similar to that of Gurobi and Mosek. Cplex and Scip emerge as the least competitive solvers.
The results provide guidelines on how each solver behaves on this class of problems and highlight the importance of choosing a solver suited to the problem type.


\section{Introduction}
\label{ch:intro}



The congested capacitated facility location problem (CCFLP) is a discrete optimization problem which holds the interest of the scientific community due to its practical applicability and the theoretical challenges it offers. Introduced in 1995 in \cite{desrochers1995congested}, it represents an challenging variant of the facility location problems family, integrating elements of congestion and capacity.
The CCFLP has a natural formulation as a mixed-integer convex quadratic program with indicator variables \cite{fischetti2016benders}. However, the use of perspective reformulation on such problems has been demonstrated providing a considerably stronger continuous relaxation than the natural problem formulation (see \cite{gunluk2011perspective, gunluk2010perspective} for details). Perspective reformulation of this class of problem belongs to the class of mixed-integer conic programming (MICP) and falls specifically into the case of mixed-integer second-order cone programming (MISOCP).

This paper presents a comparative analysis of commonly available solvers for mixed-integer nonlinear programming (MINLP) on the CCFLP. Our aim is to evaluate the performance of these solvers in terms of their efficiency and scalability. The solvers we consider in this study are some of the state-of-the-art algorithms widely adopted by the optimization community: Gurobi, Cplex, Mosek, Xpress, and Scip. All these solvers are commercial, except for Scip.
The evaluation of the performance of the solvers is based on three criteria: the runtime, the optimality gap at the root node, and the number of tree nodes explored. Additionally, we compare the solution profiles and the solution status returned by the solvers on the testbed. To present our results, we provide both numerical data and visual representations.


A large number of CCFLP instances with variable sizes, ranging from small to large, and different computational complexity are available in \cite{fischetti2016benders}. This enables us to test the solvers across numerous and varied scenarios of the same problem, assessing their performance as the problem size and complexity increase. 

The contribution of this paper is twofold. Primarily, we aim to provide a comprehensive comparison of five of the most common MINLP solvers on the CCFLP. As a by-product, we derive a practical guide for practitioners and researchers to assist in the choice of solvers, encouraging an informed selection based on the specific nature of the problem at hand. We remark that our goal is to assess the performance of the solvers on the CCFLP. While acknowledging the existence of comprehensive libraries of MICP problems with instances derived from a variety of different problems (e.g., \cite{friberg2016cblib}), evaluating solvers on such a broad spectrum is beyond the scope of this article.

The paper is organized as follows. In Section \ref{sec:CCFLP}, the MICP formulation of the CCFLP used as the testbed class of problems is presented. Section \ref{ch:solvers} provides a brief overview of the solvers and their solution strategies for MICP. Insights into the method chosen by each solver specifically on the testbed are also given. In Section \ref{ch:results}, we report the instances description. Furthermore, the numerical performance of the solvers is presented and discussed. Conclusions are given in Section \ref{ch:conclusions}.

\section{The MICP formulation of the CCFLP}
\label{sec:CCFLP}

In this section, we  briefly report the MICP formulation of the CCFLP proposed in \cite{fischetti2016benders} and used to obtain our testbed.

We are given a finite set $I$ of customers, each characterized by a demand $d_i \geq 0$ $i \in I$ to be served, and a finite set $J$ of potential facility locations, with a maximum serving capacity $s_j$,  $j \in J$. 
Opening a facility $j \in J$ has a cost $f_j \geq 0$, as well as serving a unit of demand of customer $i \in I$ by facility $j \in J$ has a cost $c_{ij} \geq 0$. 
Congestion at a facility $j \in J$ is measured by the load of $j$, namely the overall amount of demand served by facility $j$ towards all the served customers $i \in I$.
The congestion at open facilities may imply diseconomies of scale, e.g., employment of a larger number of overtime workers, usage of more costly materials, postponement of equipment maintenance schedules \cite{harkness2003facility}. Thus, congestion has costs to be minimized,
which are usually
modeled through quadratic convex functions. When congestion costs are modeled, a constraint to control the number of open facilities is usually considered to prevent the model from favoring a solution with many open facilities to keep their load at the minimum.
We denote by $p$ the exact number of facilities to open.
We also allow each facility to cover only a fraction of the total demand of each customer, whose total demand can thus be covered by different facilities.

The CCFLP consists of determining the facility locations to open and the fraction of demand  served by the selected facilities for each customer in order to minimize the overall cost obtained as the sum of facility opening, customer allocation, and congestion costs.
Let  $y_j\in\{0,1\}$ denote the closing/opening of the facility location $j\in J$ and $x_{ij} \geq 0$ for each $i \in I, j \in J$  the fraction of the demand of customer $i \in I$ served by facility $j \in J$ when $j$ is open. 
Thus the overall cost  is the sum of facility opening costs, customer allocation costs, and congestion costs, and it is written as
\cite{fischetti2016benders}
$$ \sum_{j \in J} f_{j}y_{j} + \sum_{i \in I}\sum_{j \in J} d_i c_{ij} x_{ij} + b \sum_{j \in J} v_j + a \sum_{j \in J} v_j^2$$
where  $v_j= \sum_{i \in I} d_i x_{ij} $ denote the exact demand served by facility $j \in J$, and $a$ and $b$ are non-negative input coefficients.

Starting from the natural formulation of the CCFLP as a mixed-integer quadratic convex program, we can derive its perspective reformulation in which the quadratic term in the objective function (representing the congestion) is replaced by a non-negative variable $z$ so as to move the non-linearity from the objective to the constraints. We refer the reader to \cite{fischetti2016benders} to see precisely how the perspective reformulation is derived from the original problem formulation. 
The CCFLP can be modeled using its perspective formulation as the following mixed-integer conic programming:
\begin{align}
\min_{x,y,v,z} \quad &  \sum_{j \in J} f_{j}y_{j} + \sum_{i \in I}\sum_{j \in J} d_i c_{ij} x_{ij} + b \sum_{j \in J} v_j + a \sum_{j \in J} z_j \nonumber\\
& \sum_{j \in J}y_{j} = p \label{pmedian_constraint}\\
& v_j^2 \leq z_jy_j && j \in J \label{perspective_constraints}\\
& \sum_{i \in I} d_i x_{ij} = v_j && j \in J \nonumber\\
&\sum_{j \in J} x_{ij} = 1 && i \in I \label{assignment_constraint}\\
& x_{ij} \leq y_j && i \in I, j \in J \label{VUB_constraint}\\
& v_j \leq s_j y_j && j \in J \label{capacity_constraint}\\
& x_{ij} \geq 0 && i \in I, j \in J \nonumber\\
&z_j \geq 0 && j \in J \nonumber\\
&y_{j} \in \{0,1\} && j \in J. \nonumber
\end{align}

The rotated second-order cone constraints \eqref{perspective_constraints}, together with the objective function's minimization, guarantee the quadratic load $v^2_j$ of a facility $j$ is zero if the facility is closed, and it is $z_j$ if the facility is open.
Assignment constraints \eqref{assignment_constraint} guarantee that all customer demand is met,
while capacity constraints \eqref{capacity_constraint} ensure that the capacity of an open facility is not exceeded.
Constraints \eqref{VUB_constraint} state that allocation to a facility $j$ is only possible if it is open; they are known to strengthen the formulation significantly. 
Constraint \eqref{pmedian_constraint} controls the number of open facilities, avoiding that the optimal solution requires opening too many facilities.

Given the definition of a $3$-dimensional rotated quadratic cone
$$\mathcal K_r^3(j) = \{(z_j,y_j,v_j)\in \mathbb R^3 \text{ } |\text{ }  2z_jy_j \geq v_j^2, \text{ } z_j \geq 0, \text{ }y_j \geq 0\} \quad \text{for } j  \in J,$$
we can represent constraint \eqref{perspective_constraints} equivalently as 
\begin{equation}\label{cones_mosek}
    (z_j,y_j,\sqrt{2}v_j) \in \mathcal K_r^3(j)
\end{equation}
for $j \in J$. Note that \eqref{cones_mosek} also includes the non-negativity of $z$. Mosek does not automatically detect second-order cones, and it requires the implementation of the cones as in \eqref{cones_mosek}. All the other solvers admit the implementation of quadratic cones as in \eqref{perspective_constraints}.

\section{Solvers at a glance}\label{ch:solvers}

In this section, we briefly introduce the five software for mixed-integer nonlinear programming and describe their main characteristics. Indeed, most of the solvers combine several techniques to improve
their performance. 
In particular, MICP problems can be addressed through a branch-and-bound strategy based on either a continuous linear programming (LP) solver or a continuous quadratically
constrained programming (QCP) solver. 
The former has the advantage of faster node processing but leads to weaker relaxations at the nodes; the latter yields tighter bounds at the expense of a longer node processing time (see, e.g., \cite{CplexWebsite}). The branch-and-bound strategy based on the LP solver is an outer-approximation approach (OA) in which first-order approximations of the quadratic constraints are used to replace the quadratic constraints. The linear approximations are dynamically refined throughout the process to yield an optimal solution that satisfies all the quadratic constraints of the model.

In this section, we discuss the default solution approach implemented by each solver to tackle MICPs in general 
and our testbed in particular. 
The summary below is not intended to give an in-depth
analysis, but to better exemplify the similarities and differences between the solvers. We refer the reader to the solver website for more details. The information reported has been extracted from the available documentation of each solver, or in case documentation was not complete, from the log file produced on our instances, as it can be considered a valid statement on this specific testbed.
In particular, we aim to enlighten the node strategy (LP or QCP) used by the solver. 
  We also report the algorithm used at the root node and how the root node relaxation is obtained. More specifically, at the root node each solver involves multiple procedures to improve the bounds, such as generating cuts and applying heuristics; consequently, each solver reports a sequence of improving values of the bounds at the root node. In the description below, we discuss the very first relaxation at the root node and the corresponding bound that can be derived from the log file. The value of the initial lower bound at the root node is mainly improved by cut generation. Details of the approaches to improve this bound are not recoverable from the log file.

The solvers can be used with different interfaces. We set the default setting for all the solvers and used their Python interface.
 For each solver, we provide the type of license available, the version and the python package we used in the numerical experiments, and the URL to the solver website. 

 Table \ref{tab:solvers} provides a summary of the solver features.

\subsection{Gurobi}

\textbf{Name:} Gurobi
\\\textbf{License type:} Free for academic use; commercial
\\\textbf{Version:} 10.0.0
\\\textbf{Python package:} gurobipy
\\\textbf{Website:} \url{https://www.gurobi.com/}

\medskip 
\noindent
In the default setting, Gurobi chooses automatically between a linearized OA approach and a branch-and-bound strategy based on continuous QCP relaxations at each node \cite{gurobiWebsite}. Its strategy can be dynamically adjusted during the optimization. 

We can  infer from the log files that the very initial root bound is obtained by relaxing both integrality and nonlinear constraints. Further,  Gurobi uses the dual simplex to solve the very initial root relaxation. 

\subsection{Cplex}

\textbf{Name:} IBM ILOG CPLEX Optimization Studio
\\\textbf{License type:} Free for academic use; commercial
\\\textbf{Version:} 20.1.0
\\\textbf{Python package:} docplex
\\\textbf{Website:} \url{https://www.ibm.com/products/ilog-cplex-optimization-studio}

\medskip 
\noindent
According to the available documentation \cite{CplexWebsite}, Cplex can exploit both the LP-based and the QCP-based strategies to solve MIQCP formulations. 

We can  infer from the log files that the very initial root bound is obtained by relaxing both integrality and nonlinear constraints.
After that,  Cplex applies OA at the root node in all instances we tested. The linear approximation cuts are denoted as cone linearizations in the Cplex log.

\subsection{Mosek}

\textbf{Name:} MOSEK
\\\textbf{Version:} 9.3.20
\\\textbf{License type:} Free for academic use; commercial
\\\textbf{Python package}: mosek
\\\textbf{Website:} \url{https://www.mosek.com/}

\medskip 
\noindent
Mosek requires the implementation of the conic constraints as in \eqref{cones_mosek}. In the default setting, it implements an efficient interior point (IP) algorithm for conic problems described in \cite{mosekWebsite}.
OA can be applied when solving relaxations of conic problems only if activated through an option.

According to the tests we run, the value of the very initial root relaxation corresponds to relaxing only the integrality constraints.

\subsection{Xpress}

\textbf{Name:} FICO XPRESS
\\\textbf{Version:} 8.13.4
\\\textbf{License type:} Free for academic use; commercial
\\\textbf{Python package:} xpress
\\\textbf{Website:} \url{https://www.fico.com/en/products/fico-xpress-optimization}

\medskip 
\noindent

Xpress documentation \cite{xpressWebsite} declares that mixed-integer quadratically constrained and mixed-integer second-order cone problems are solved via a traditional branch-and-bound using either the barrier method to solve the node problems or via OA. In the default setting, the method used is automatically determined. It is also stated that it may be beneficial to solve the root node of such problems by the barrier, even if OA is used in the branching nodes. 

On our testbed, by default, the Newton barrier method is used to solve the very initial root relaxation of the model, and then OA is applied (at least at the root node). The lower bound obtained with the barrier corresponds to relaxing only integrality constraints.
After the barrier, Xpress performs a crossover phase to obtain a basic optimal solution with the simplex method, as the solution of the barrier may not lie exactly on the boundary of the feasible region \cite{xpress_crossover}. The bound found after the crossover is used as starting value of the objective relaxation in the tree search. The crossover deteriorates the value of the bound found by the barrier in most of the instances we tested. Nevertheless, the bound after the crossover is always better than the one obtained by relaxing both the integrality and the nonlinear constraints. In this sense, the combination of barrier and simplex methods proves to be effective.

\subsection{Scip}

\textbf{Name:} SCIP
\\\textbf{Version:} 8.0.1
\\\textbf{License type:} Free for academic use
\\\textbf{Python package:} pyscipopt
\\\textbf{Website:} \url{https://www.scipopt.org/}

\medskip 
\noindent
Scip is the only non-commercial solver among the ones we tested. 
Scip implements a spatial branch-and-bound algorithm based on a linear OA \cite{vigerske2018Scip}, i.e., it applies only LP methods at the nodes.

On our testbed, the value of the very initial root relaxation corresponds to the value obtained by relaxing both integrality and nonlinear constraints.

\begin{table}[ht]
\centering
\caption{Comparison over solvers: T indicates what is done on our testbed; \protect\tikzcmark \hspace{0.01cm} indicates what the solver can do according to the documentation.}
\begin{tabular}{lcccccc} \toprule
& \multicolumn{2}{c}{\textbf{ Constraint Relaxation at Root}} & \multicolumn{2}{c}{\textbf{Node Solver}}        &             \\
\textbf{Solver} &  Integrality Constraint &  Quadratic Constraint & LP & QCP & \\ \cmidrule(lr){1-1}\cmidrule(lr){2-3}\cmidrule(lr){4-5}
Gurobi & T & T & \tikzcmark & \tikzcmark \\
Cplex & T & T & T & \tikzcmark \\
Mosek & T & & \tikzcmark & T \\
Xpress &T & & T & \tikzcmark\\
Scip & T & T & T & \\
\bottomrule
\end{tabular} \label{tab:solvers}
\end{table}

\section{Numerical Experiments}
\label{ch:results}
In this section, we give details on the CCFLP instances used as a testbed and analyze the performance of the five solvers. All computations were performed on a Ubuntu server with an Intel(R) Xeon(R) Gold 5218 CPU with 96 GB RAM.

For all the solvers, the termination criteria were a time limit of 14\,400 s (4 hours), or an optimality gap below $10^{-6}$. The following tolerances were also set: optimality tolerance below $10^{-6}$, integrality tolerance below $10^{-6}$, and feasibility tolerance below $10^{-4}$. For all the other options, the default setting was used in each solver. 

Like in previous benchmark studies \cite{dolan2002benchmarking,neumaier2005comparison}, we recognize that testing optimization software is a difficult problem and that there may be objections to the testing presented in this report. Indeed, the performance of a solver may improve significantly if non-default options are given.
However, the definition of a customized setting for each solver is extremely time-consuming and itself represents a difficult problem to pursue. Such customization requires deep knowledge of both the problem under study and the specific solver, and it might not be viable. Thus, although the conclusions  could be biased on our choice of using the default setting, we guess that  it should not affect comparisons that rely on large differences.

\subsection{Instances}

The instances used as testbed have been taken from those used in \cite{fischetti2016benders} and have been kindly provided by one of its authors. The original testbed is made of 60 instances of the CCFLP with up to 1000 facilities and 1000 customers. We selected 30 instances out of the 60 available as the testing process turned out to be extremely time-consuming due to the fact we needed to run each instance on five solvers, and each instance could take up to four hours (which is the time limit we set) for its solution.

The instances we used can be classified into five blocks, each consisting of six instances of the same size $|I|$, $|J|$. Each block, in turn, is divided into two sub-blocks which differ in the ratio $ r = \sum_{j \in J} s_j/\sum_{i \in I} d_i$ between the overall serving capacity and the demand. 
Each sub-block is composed of three instances that differ from each other in  the number $p$ of open facilities required, which is obtained as a fraction $\pi$ of the number of available facilities $|J|$ ($p = \lfloor \pi |J| \rfloor$). Table \ref{tab:instances} reports the characteristics of the instances. The ID number reported in Table \ref{tab:instances} will be used in Section \ref{ch_results} to identify the specific instance. 
The instances are sorted in Table \ref{tab:instances} by increasing size, with the exception of the last block, which contains medium-sized but more complex problems as shown by the computational results in Section \ref{ch_results} and in \cite{fischetti2016benders}. For further description of instances, see \cite{fischetti2016benders}.

\begin{table}[h]
\caption{Testbed features. For each instance under column \textbf{ID}, the identifier \textbf{Inst.} used in \cite{fischetti2016benders}, the number of customers \pmb{$|I|$}, the number of facilities \pmb{$|J|$}, the scaling factor \pmb{$r$}, and the parameter \pmb{$\pi$} are given. }
\label{tab:instances}
\begin{center}
\resizebox{0.75\textwidth}{!}{
\begin{tabular}{lrrrrrrlrrrrr}
\toprule
\textbf{ID}  &		 \textbf{Inst.} &	 \pmb{$|I|$} &	 \pmb{$|J|$} &	 \pmb{$r$} &	 \pmb{$\pi$} & & \textbf{ID}  &	 \textbf{Inst.} &	 \pmb{$|I|$} &	 \pmb{$|J|$} &	 \pmb{$r$} & \pmb{$\pi$} \\
\cmidrule(lr){1-6}\cmidrule(lr){8-13}
1	&	1 & 300 & 300 & 5 & 0,4	&	&		16	&	6 & 700 & 700 & 10 & 0,4			\\
2	&	1 & 300 & 300 & 5 & 0,6	&	&		17	&	6 & 700 & 700 & 10 & 0,6			\\
3	&	1 & 300 & 300 & 5 & 0,8	&	&		18	&	6 & 700 & 700 & 10 & 0,8			\\ \cmidrule(lr){8-13}
4	&	6 & 300 & 300 & 10 & 0,4	&	&		19	&	1 & 1000 & 1000 & 5 & 0,4			\\
5	&	6 & 300 & 300 & 10 & 0,6	&	&		20	&	1 & 1000 & 1000 & 5 & 0,6			\\
6	&	6 & 300 & 300 & 10 & 0,8	&	&		21	&	1 & 1000 & 1000 & 5 & 0,8			\\  \cmidrule(lr){1-7}
7	&	1 & 500 & 500 & 5 & 0,4	&	&		22	&	11 & 1000 & 1000 & 15 & 0,4			\\
8	&	1 & 500 & 500 & 5 & 0,6	&	&		23	&	11 & 1000 & 1000 & 15 & 0,6			\\
9	&	1 & 500 & 500 & 5 & 0,8	&	&		24	&	11 & 1000 & 1000 & 15 & 0,8 			\\\cmidrule(lr){8-13}
10	&	11 & 500 & 500 & 15 &  0,4	&	&		25	&	1 & 300 & 1500 & 5 & 0,4			\\
11	&	11 & 500 & 500 & 15 & 0,6	&	&		26	&	1 & 300 & 1500 & 5 & 0,6			\\
12	&	11 & 500 & 500 & 15 & 0,8	&	&		27	&	1 & 300 & 1500 & 5 & 0,8			\\ \cmidrule(lr){1-7}
13	&	1 & 700 & 700 & 5 &  0,4	&	&		28	&	16 & 300 & 1500 & 20 & 0,4			\\
14	&	1 & 700 & 700 & 5 & 0,6	&	&		29	&	16 & 300 & 1500 & 20 & 0,6			\\
15	&	1 & 700 & 700 & 5 & 0,8	&	&		30	&	16 & 300 & 1500 & 20 & 0,8	\\\bottomrule
\end{tabular}}
\end{center}
\end{table}
 
\subsection{Results}\label{ch_results}

The purpose of this section is to compare the performance of the solvers on the MICP formulation of the CCFLP. The metrics we use for the comparison are: i) the runtime, i.e. the  computational time expressed in seconds, ii) the relative gap at the end of root node analysis,
iii) the relative gap at the end of the optimization
and iv) the number of branching nodes explored. 
The results are collected in Table \ref{tab:results}.

\begin{figure}[h]
    \centering
    \includegraphics[width=0.9\linewidth]{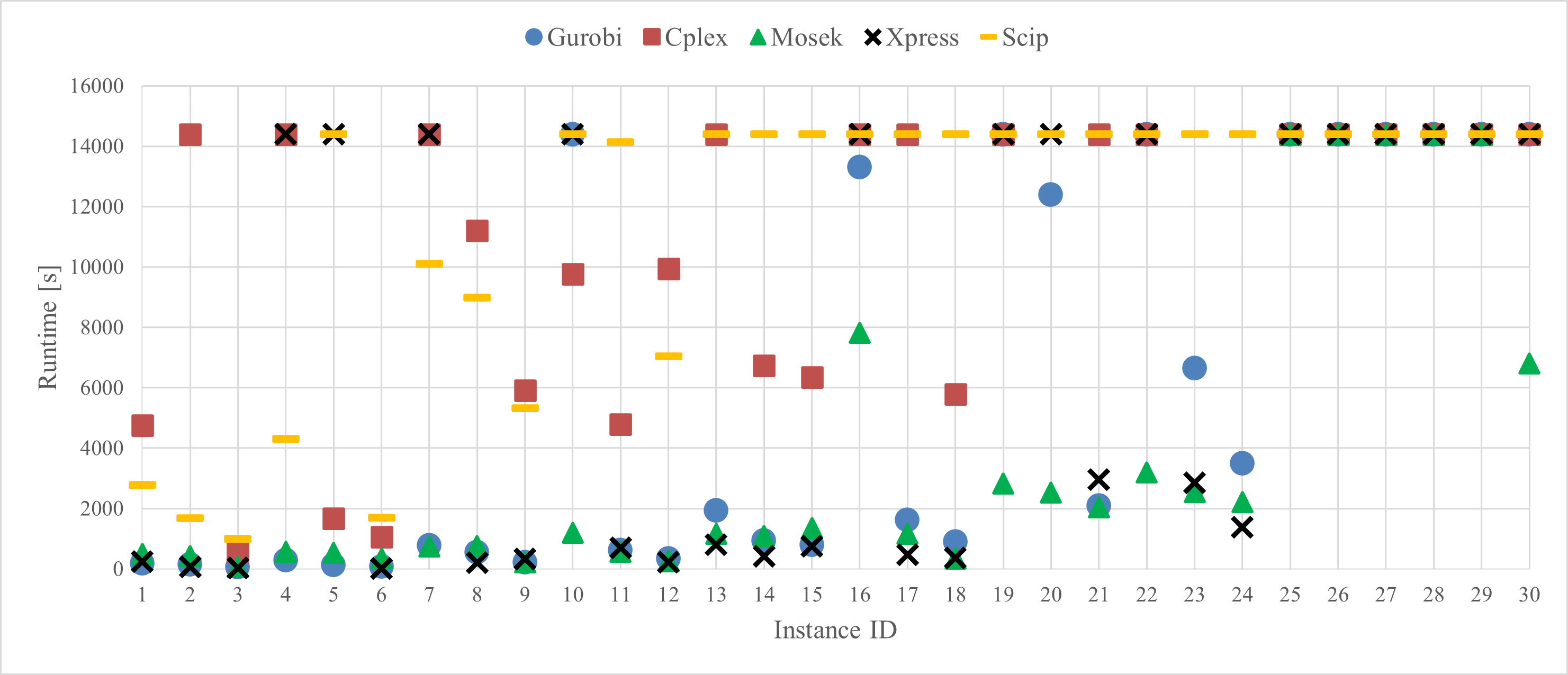}
    \caption{Runtime plot.}
    \label{fig:runtime}
\end{figure}

\begin{figure}[h]
    \centering
    \includegraphics[width=0.9\linewidth]{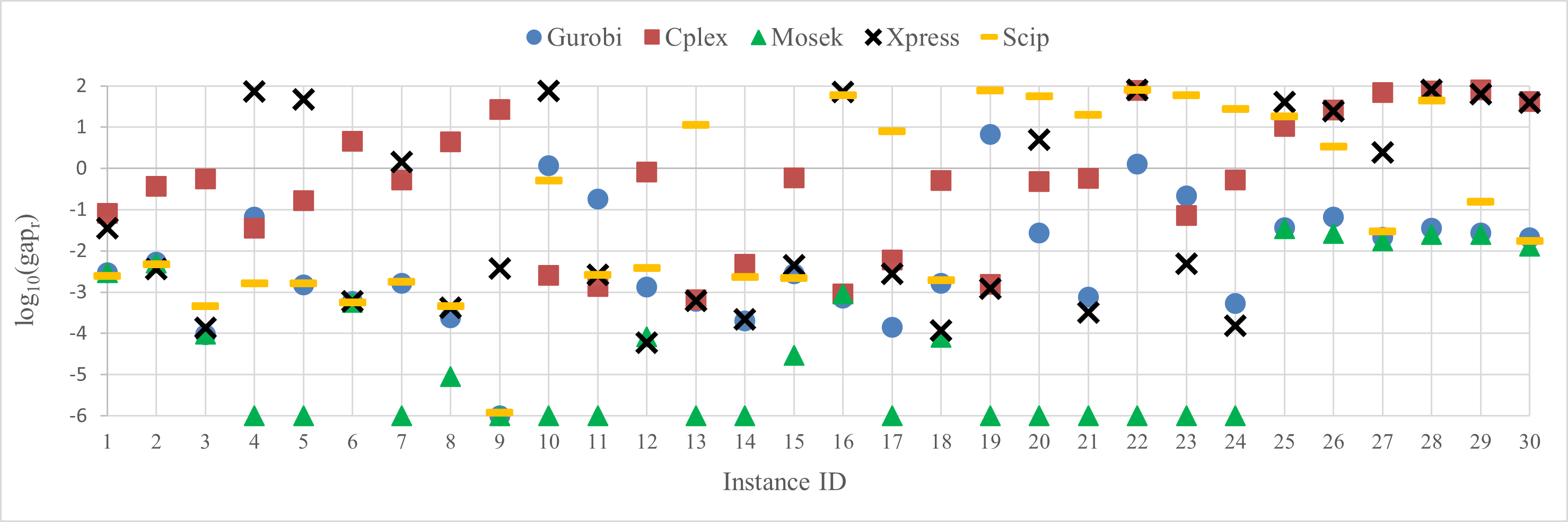}
    \caption{Logarithmic relative root gap plot.}
    \label{fig:rootgap}
\end{figure}

\begin{figure}[h]
    \centering
    \includegraphics[width=0.9\linewidth]{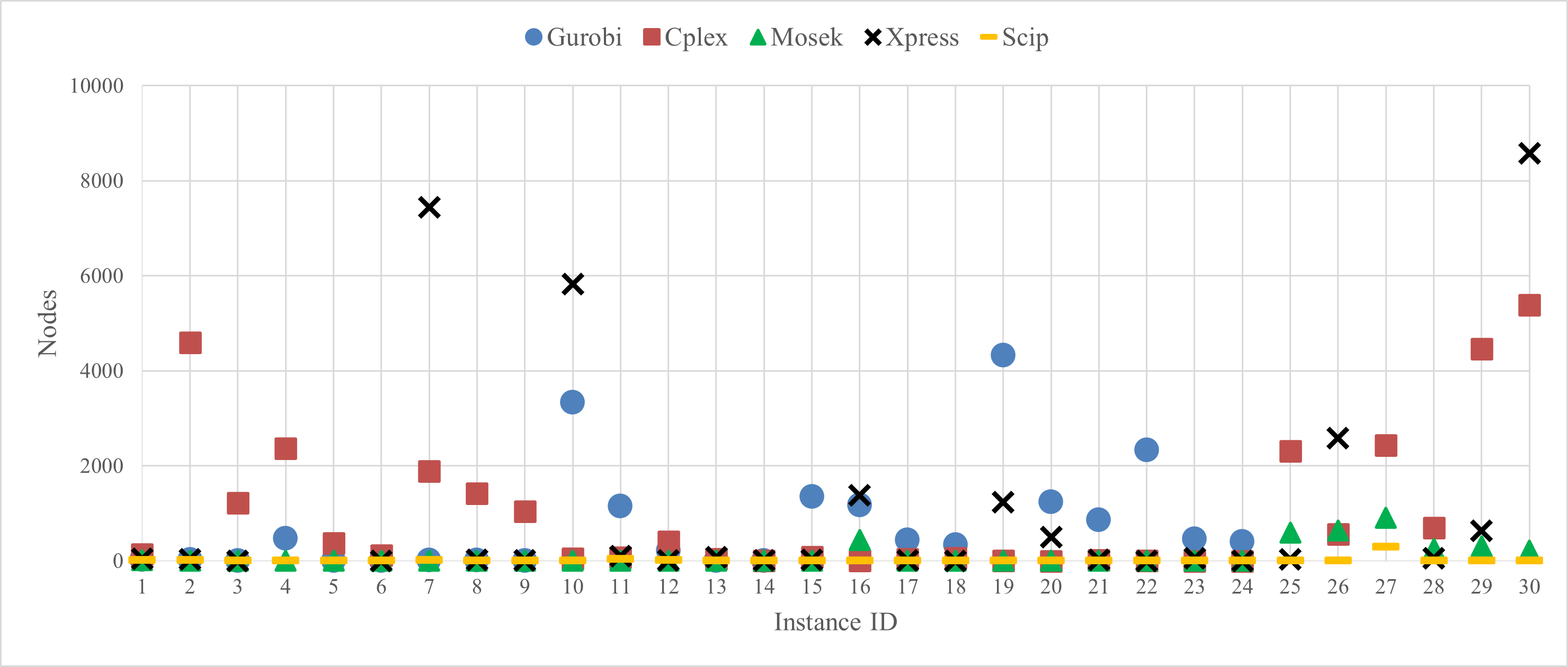}
    \caption{Number of explored nodes plot.}
    \label{fig:nodes}
\end{figure}

As for the comparison based on metrics  i), ii) and iv), we present two types of plots for each metric based on the results of the instances.
In particular, for each of the three metrics, Figures \ref{fig:runtime}, \ref{fig:rootgap} and \ref{fig:nodes},
are two-dimensional plots, with the instances on the x-axis and the value of the chosen metric on the y-axis for each solver. 
In the Figures \ref{fig:runtime}, \ref{fig:rootgap} and \ref{fig:nodes}
each solver is represented by a symbol
and color. 
This type of plot, displaying the solvers' value of the selected metric across all instances, is useful in assessing how the solver's performance is impacted by the size/complexity of the instance.

\begin{table}[t!]
\setlength{\belowcaptionskip}{-20pt}
\caption{Solvers performance on the testbed. Column \textbf{gap\textsubscript{r\textsubscript{i}}} gives the initial relative gap at the root node, \textbf{gap\textsubscript{r}} the relative gap at the end of root node analysis, \textbf{gap} the relative gap at the end of the optimization, \textbf{runtime} the computational time expressed in seconds, \textbf{nodes} the number of branching nodes. Gaps are evaluated with respect to the best upper bound. Runs reaching the time limit of 14\,400 seconds (4 hours) are indicated by “TL”; runs aborted due to memory limit are indicated by “*.”}
\begin{center}
\resizebox{\textwidth}{!}{
\begin{tabular}{lrrrrrrrrrrrrrrr}
\toprule
     & \multicolumn{3}{c}{\textbf{Gurobi}}      & \multicolumn{3}{c}{\textbf{Cplex}}   
     & \multicolumn{3}{c}{\textbf{Mosek}}       
     & \multicolumn{3}{c}{\textbf{Xpress}}      & \multicolumn{3}{c}{\textbf{Scip}}\\   \cmidrule(lr){2-4}\cmidrule(lr){5-7}\cmidrule(lr){8-10}\cmidrule(lr){11-13}\cmidrule(lr){14-16}
\textbf{ID} & \textbf{gap\textsubscript{r\textsubscript{i}}{[}\%{]}} & \textbf{gap\textsubscript{r}{[}\%{]}} & \textbf{gap{[}\%{]}} & \textbf{gap\textsubscript{r\textsubscript{i}}{[}\%{]}} & \textbf{gap\textsubscript{r}{[}\%{]}} & \textbf{gap{[}\%{]}} &  \textbf{gap\textsubscript{r\textsubscript{i}}{[}\%{]}} & \textbf{gap\textsubscript{r}{[}\%{]}} & \textbf{gap{[}\%{]}} &  \textbf{gap\textsubscript{r\textsubscript{i}}{[}\%{]}} & \textbf{gap\textsubscript{r}{[}\%{]}} & \textbf{gap{[}\%{]}} & \textbf{gap\textsubscript{r\textsubscript{i}}{[}\%{]}} & \textbf{gap\textsubscript{r}{[}\%{]}} & \textbf{gap{[}\%{]}}\\
\cmidrule(lr){1-1}\cmidrule(lr){2-4}\cmidrule(lr){5-7}\cmidrule(lr){8-10}\cmidrule(lr){11-13}\cmidrule(lr){14-16}
1  & 8,50E+01 & 2,95E-03 & 5,16E-05 & 8,50E+01 & 8,27E-02 & 0,00E+00 & 3,00E-03 & 3,00E-03 & 0,00E+00 & 7,90E+01 & 3,59E-02 & 2,11E-05 & 8,50E+01 & 2,45E-03 & 5,89E-05 \\
2  & 6,78E+01 & 5,37E-03 & 0,00E+00 & 6,78E+01 & 3,66E-01 & 1,21E-01 & 5,00E-03 & 5,00E-03 & 0,00E+00 & 6,83E-01 & 3,61E-03 & 9,93E-05 & 6,78E+01 & 4,74E-03 & 0,00E+00 \\
3  & 4,83E+01 & 9,34E-05 & 9,34E-05 & 4,83E+01 & 5,63E-01 & 4,38E-05 & 9,53E-05 & 9,53E-05 & 9,53E-05 & 1,35E+01 & 1,39E-04 & 8,65E-05 & 4,83E+01 & 4,47E-04 & 0,00E+00 \\
4  & 8,48E+01 & 6,62E-02 & 6,59E-05 & 8,48E+01 & 3,61E-02 & 9,93E-03 & 4,41E-03 & 0,00E+00 & 0,00E+00 & 8,39E+01 & 7,26E+01 & 9,73E+00 & 8,48E+01 & 1,60E-03 & 0,00E+00 \\
5  & 6,81E+01 & 1,52E-03 & 8,24E-05 & 6,81E+01 & 1,66E-01 & 0,00E+00 & 1,55E-03 & 0,00E+00 & 0,00E+00 & 6,69E+01 & 4,66E+01 & 1,00E+01 & 6,81E+01 & 1,60E-03 & 7,96E-05 \\
6  & 5,06E+01 & 5,98E-04 & 0,00E+00 & 5,06E+01 & 4,62E+00 & 0,00E+00 & 5,73E-04 & 5,73E-04 & 0,00E+00 & 5,94E-04 & 5,93E-04 & 1,00E-04 & 5,06E+01 & 5,56E-04 & 0,00E+00 \\
7  & 8,48E+01 & 1,63E-03 & 6,88E-05 & 8,48E+01 & 5,24E-01 & 5,22E-01 & 2,60E-03 & 0,00E+00 & 0,00E+00 & 8,20E+01 & 1,43E+00 & 1,36E-01 & 8,48E+01 & 1,76E-03 & 0,00E+00 \\
8  & 6,80E+01 & 2,38E-04 & 0,00E+00 & 6,80E+01 & 4,39E+00 & 9,77E-05 & 6,74E-04 & 9,04E-06 & 9,04E-06 & 2,13E+00 & 4,27E-04 & 9,79E-05 & 6,80E+01 & 4,60E-04 & 0,00E+00 \\
9  & 4,91E+01 & 0,00E+00 & 0,00E+00 & 4,91E+01 & 2,69E+01 & 2,55E-05 & 0,00E+00 & 0,00E+00 & 0,00E+00 & 3,82E+01 & 3,71E-03 & 6,18E-05 & 4,91E+01 & 1,20E-06 & 1,20E-06 \\
10 & 8,45E+01 & 1,16E+00 & 7,11E-02 & 8,45E+01 & 2,58E-03 & 1,48E-05 & 2,51E-03 & 0,00E+00 & 0,00E+00 & 8,38E+01 & 7,49E+01 & 6,66E+01 & 8,45E+01 & 5,00E-01 & 5,00E-01 \\
11 & 6,74E+01 & 1,82E-01 & 0,00E+00 & 6,74E+01 & 1,38E-03 & 5,61E-05 & 1,36E-03 & 0,00E+00 & 0,00E+00 & 3,07E+01 & 2,62E-03 & 6,93E-05 & 6,74E+01 & 2,56E-03 & 0,00E+00 \\
12 & 4,88E+01 & 1,33E-03 & 0,00E+00 & 4,88E+01 & 8,23E-01 & 3,43E-06 & 8,33E-05 & 8,33E-05 & 8,33E-05 & 2,91E+00 & 5,98E-05 & 2,41E-05 & 4,88E+01 & 3,79E-03 & 0,00E+00 \\
13 & 8,40E+01 & 5,92E-04 & 4,66E-05 & 8,40E+01 & 6,51E-04 & 5,44E-04 & 6,43E-04 & 0,00E+00 & 0,00E+00 & 8,35E+00 & 6,17E-04 & 9,17E-05 & 8,40E+01 & 1,12E+01 & 1,12E+01 \\
14 & 6,61E+01 & 2,03E-04 & 9,54E-05 & 6,61E+01 & 4,83E-03 & 6,00E-05 & 2,48E-04 & 0,00E+00 & 0,00E+00 & 4,28E+00 & 2,20E-04 & 9,57E-05 & 6,61E+01 & 2,31E-03 & 2,31E-03 \\
15 & 4,71E+01 & 2,79E-03 & 0,00E+00 & 4,71E+01 & 5,86E-01 & 9,89E-05 & 5,11E-04 & 2,91E-05 & 2,91E-05 & 2,16E+01 & 4,50E-03 & 7,84E-05 & 4,71E+01 & 2,20E-03 & 2,20E-03 \\
16 & 8,51E+01 & 7,17E-04 & 0,00E+00 & 8,51E+01 & 9,00E-04 & 7,16E-04 & 9,10E-04 & 9,10E-04 & 0,00E+00 & 8,41E+01 & 7,12E+01 & 6,03E+01 & 8,51E+01 & 5,93E+01 & 5,93E+01 \\
17 & 6,75E+01 & 1,43E-04 & 3,90E-05 & 6,75E+01 & 6,07E-03 & 9,80E-04 & 1,36E-04 & 0,00E+00 & 0,00E+00 & 1,66E+01 & 2,82E-03 & 8,34E-05 & 6,75E+01 & 7,84E+00 & 7,84E+00 \\
18 & 4,87E+01 & 1,63E-03 & 0,00E+00 & 4,87E+01 & 5,07E-01 & 4,69E-05 & 7,96E-05 & 7,96E-05 & 7,96E-05 & 1,63E+01 & 1,20E-04 & 9,28E-05 & 4,87E+01 & 1,92E-03 & 1,92E-03 \\
19 & 8,41E+01 & 6,74E+00 & 6,73E+00 & 8,41E+01 & 1,55E-03 & 1,48E-03 & 1,34E-03 & 0,00E+00 & 0,00E+00 & 1,41E+01 & 1,21E-03 & 7,11E-04 & 8,41E+01 & 7,76E+01 & 7,76E+01 \\
20 & 6,63E+01 & 2,76E-02 & 0,00E+00 & 6,63E+01 & 4,75E-01 & *        & 3,14E-04 & 0,00E+00 & 0,00E+00 & 6,39E+01 & 4,97E+00 & 1,16E+00 & 6,63E+01 & 5,52E+01 & 5,52E+01 \\
21 & 4,74E+01 & 7,68E-04 & 0,00E+00 & 4,74E+01 & 5,68E-01 & 5,67E-01 & 3,22E-04 & 0,00E+00 & 0,00E+00 & 6,99E+00 & 3,27E-04 & 9,30E-05 & 4,74E+01 & 1,99E+01 & 1,99E+01 \\
22 & 8,32E+01 & 1,29E+00 & 4,44E-04 & 8,32E+01 & 7,90E+01 & 7,90E+01 & 9,73E-04 & 0,00E+00 & 0,00E+00 & 8,07E+01 & 8,07E+01 & 8,07E+01 & 8,32E+01 & 7,95E+01 & 7,95E+01 \\
23 & 6,48E+01 & 2,19E-01 & 5,08E-05 & 6,48E+01 & 7,27E-02 & *        & 2,17E-04 & 0,00E+00 & 0,00E+00 & 2,91E+01 & 4,86E-03 & 7,09E-05 & 6,48E+01 & 5,85E+01 & 5,85E+01 \\
24 & 4,62E+01 & 5,36E-04 & 0,00E+00 & 4,62E+01 & 5,30E-01 & *        & 1,93E-04 & 0,00E+00 & 0,00E+00 & 3,51E+00 & 1,53E-04 & 9,65E-05 & 4,62E+01 & 2,73E+01 & 2,73E+01 \\
25 & 9,67E+01 & 3,63E-02 & 2,14E-02 & 9,67E+01 & 1,05E+01 & 1,05E+01 & 3,43E-02 & 3,43E-02 & 3,01E-02 & 9,54E+01 & 4,08E+01 & 1,83E+01 & 9,67E+01 & 1,79E+01 & 1,79E+01 \\
26 & 9,48E+01 & 6,60E-02 & 5,53E-02 & 9,48E+01 & 2,61E+01 & 2,61E+01 & 2,65E-02 & 2,64E-02 & 2,11E-02 & 9,32E+01 & 2,42E+01 & 4,09E+00 & 9,48E+01 & 3,37E+00 & 3,37E+00 \\
27 & 9,22E+01 & 2,12E-02 & 8,75E-03 & 9,22E+01 & 6,88E+01 & 2,11E+01 & 1,75E-02 & 1,75E-02 & 1,04E-02 & 8,77E+01 & 2,47E+00 & 1,28E-01 & 9,22E+01 & 2,98E-02 & 2,98E-02 \\
28 & 9,68E+01 & 3,52E-02 & 2,20E-02 & 9,68E+01 & 7,62E+01 & 4,88E+01 & 2,46E-02 & 2,46E-02 & 1,78E-02 & 9,63E+01 & 8,11E+01 & 6,32E+01 & 9,68E+01 & 4,43E+01 & 4,43E+01 \\
29 & 9,49E+01 & 2,75E-02 & 1,15E-02 & 9,49E+01 & 7,96E+01 & 6,16E+01 & 2,51E-02 & 2,51E-02 & 1,68E-02 & 9,34E+01 & 6,25E+01 & 3,16E+01 & 9,49E+01 & 1,53E-01 & 1,53E-01 \\
30 & 9,23E+01 & 2,11E-02 & 6,54E-03 & 9,23E+01 & 4,16E+01 & 1,60E+01 & 1,32E-02 & 1,32E-02 & 0,00E+00 & 8,73E+01 & 3,90E+01 & 1,02E+01 & 9,23E+01 & 1,72E-02 & 1,72E-02\\
 \bottomrule\\\\
\end{tabular}}

\resizebox{0.8\textwidth}{!}{
\begin{tabular}{lrrrrrrrrrr}
\toprule
     & \multicolumn{2}{c}{\textbf{Gurobi}}      & \multicolumn{2}{c}{\textbf{Cplex}}   
     & \multicolumn{2}{c}{\textbf{Mosek}}       
     & \multicolumn{2}{c}{\textbf{Xpress}}      & \multicolumn{2}{c}{\textbf{Scip}}\\   \cmidrule(lr){2-3}\cmidrule(lr){4-5}\cmidrule(lr){6-7}\cmidrule(lr){8-9}\cmidrule(lr){10-11}
\textbf{ID} & \textbf{runtime{[}s{]}} & \textbf{nodes} & \textbf{runtime{[}s{]}} & \textbf{nodes} &  \textbf{runtime{[}s{]}} & \textbf{nodes} & \textbf{runtime{[}s{]}} & \textbf{nodes} & \textbf{runtime{[}s{]}} & \textbf{nodes}\\
\cmidrule(lr){1-1}\cmidrule(lr){2-3}\cmidrule(lr){4-5}\cmidrule(lr){6-7}\cmidrule(lr){8-9}\cmidrule(lr){10-11}
1  & 180,68   & 29     & 4751,28   & 144  & 503,25  & 13  & 250,12  & 41    & 2767,01  & 13  \\
2  & 137,91   & 31     & TL        & 4596 & 426,85  & 10  & 69,08   & 33    & 1679,55  & 15  \\
3  & 59,25    & 1      & 584,04    & 1215 & 62,10    & 1   & 34,13   & 3     & 993,77   & 3   \\
4  & 290,06   & 472    & TL        & 2365 & 569,70  & 8   & TL      & 56354 & 4294,76  & 7   \\
5  & 120,89   & 5      & 1665,97 & 377  & 538,08  & 1   & TL      & 16377 & TL       & 3   \\
6  & 77,62    & 5      & 1049,97 & 120  & 347,20   & 6   & 21,44   & 5     & 1679,92  & 3   \\
7  & 785,39   & 19     & TL        & 1886 & 744,87  & 1   & TL      & 7438  & 10094,07 & 17  \\
8  & 556,67   & 21     & 11210,14 & 1422 & 760,33  & 6   & 200,01  & 17    & 8970,93  & 5   \\
9  & 208,44   & 1      & 5912,45 & 1044 & 212,41  & 1   & 334,78  & 11    & 5318,94  & 1   \\
10 & TL       & 3332   & 9770,40  & 49   & 1189,84 & 1   & TL      & 5825  & TL       & 1   \\
11 & 628,41   & 1148   & 4787,02   & 71   & 575,53  & 1   & 693,36  & 95    & 14124,13 & 35  \\
12 & 345,10   & 221    & 9940,87 & 411  & 246,69  & 1   & 228,00  & 12    & 7030,37  & 19  \\
13 & 1940,44  & 11     & TL        & 42   & 1168,93 & 1   & 807,00  & 75    & TL       & 1   \\
14 & 940,66   & 3      & 6728,41 & 2    & 1095,48 & 1   & 426,35  & 9     & TL       & 1   \\
15 & 789,95   & 1349   & 6341,73 & 88   & 1364,31 & 4   & 749,67  & 23    & TL       & 1   \\
16 & 13297,57 & 1176   & TL        & 5    & 7819,31 & 435 & TL      & 1385  & TL       & 1   \\
17 & 1623,74  & 436    & TL        & 43   & 1170,85 & 1   & 474,28  & 13    & TL       & 1   \\
18 & 891,67   & 341    & 5793,31 & 60   & 334,89  & 1   & 365,02  & 0     & TL       & 1   \\
19 & TL       & 4326   & TL        & 6    & 2829,77 & 1   & TL      & 1241  & TL       & 1   \\
20 & 12385,67 & 1244   & *         & *    & 2533,80  & 1   & TL      & 506   & TL       & 1   \\
21 & 2092,13  & 861    & TL        & 19   & 2040,50  & 1   & 2961,45 & 21    & TL       & 1   \\
22 & TL       & 2333   & TL        & 1    & 3194,21 & 1   & TL      & 0     & TL       & 1   \\
23 & 6638,76  & 461    & *         & *    & 2543,86 & 1   & 2862,00 & 53    & TL       & 1   \\
24 & 3488,01  & 405    & *         & *    & 2214,57 & 1   & 1376,18 & 5     & TL       & 1   \\
25 & TL       & 82392  & TL        & 2310 & TL      & 595 & TL      & 29    & TL       & 1   \\
26 & TL       & 34698  & TL        & 562  & TL      & 639 & TL      & 2588  & TL       & 1   \\
27 & TL       & 29558  & TL        & 2439 & TL      & 904 & TL      & 12430 & TL       & 298 \\
28 & TL       & 46340  & TL        & 700  & TL      & 250 & TL      & 54    & TL       & 1   \\
29 & TL       & 159342 & TL        & 4465 & TL      & 312 & TL      & 630   & TL       & 1   \\
30 & TL       & 40325  & TL        & 5386 & 6802,01 & 214 & TL      & 8575  & TL       & 10\\
\bottomrule
\end{tabular}}
\label{tab:results}
\end{center}
\end{table}

\clearpage

\begin{figure}[h]
    \centering
    \includegraphics[width=0.45\linewidth]{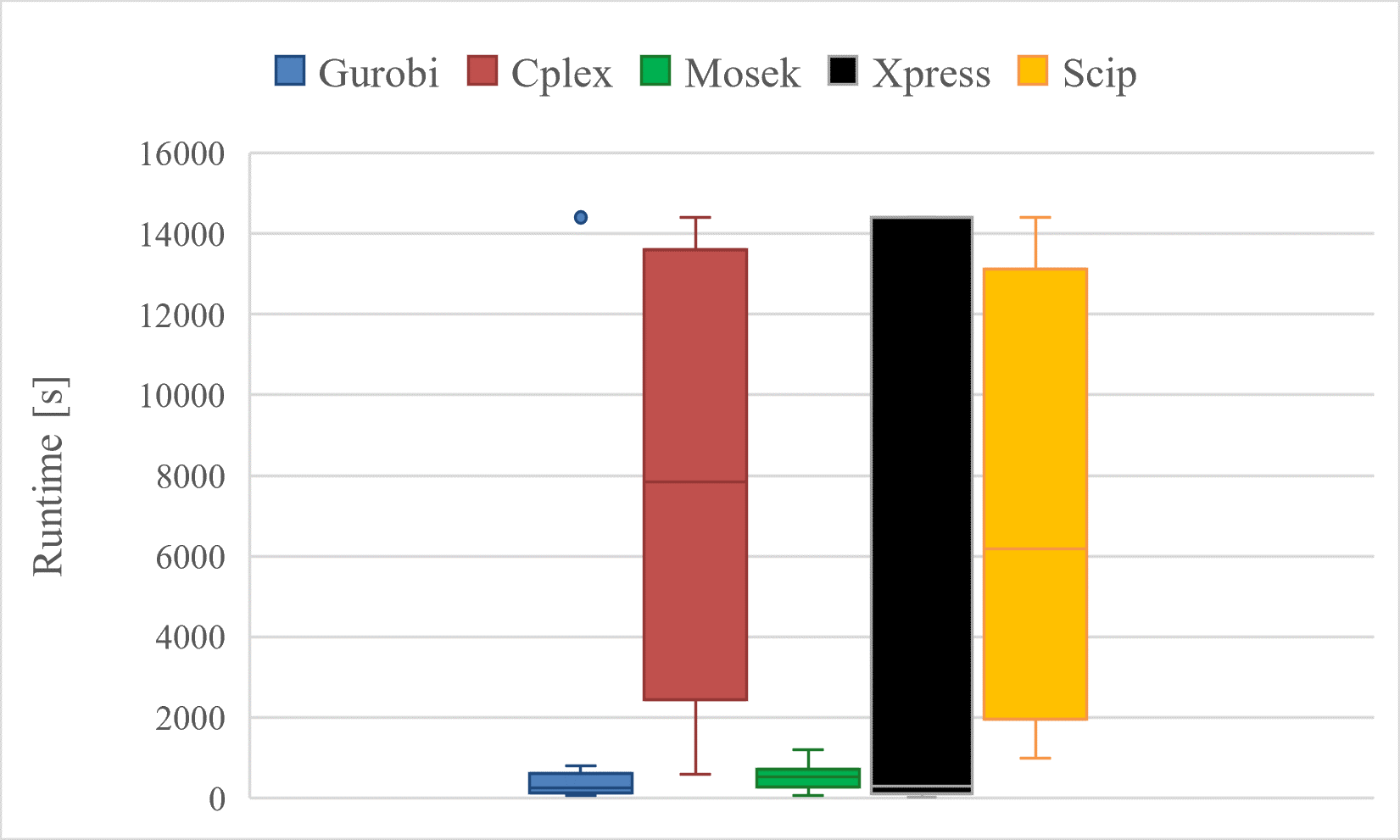}
    \includegraphics[width=0.45\linewidth]{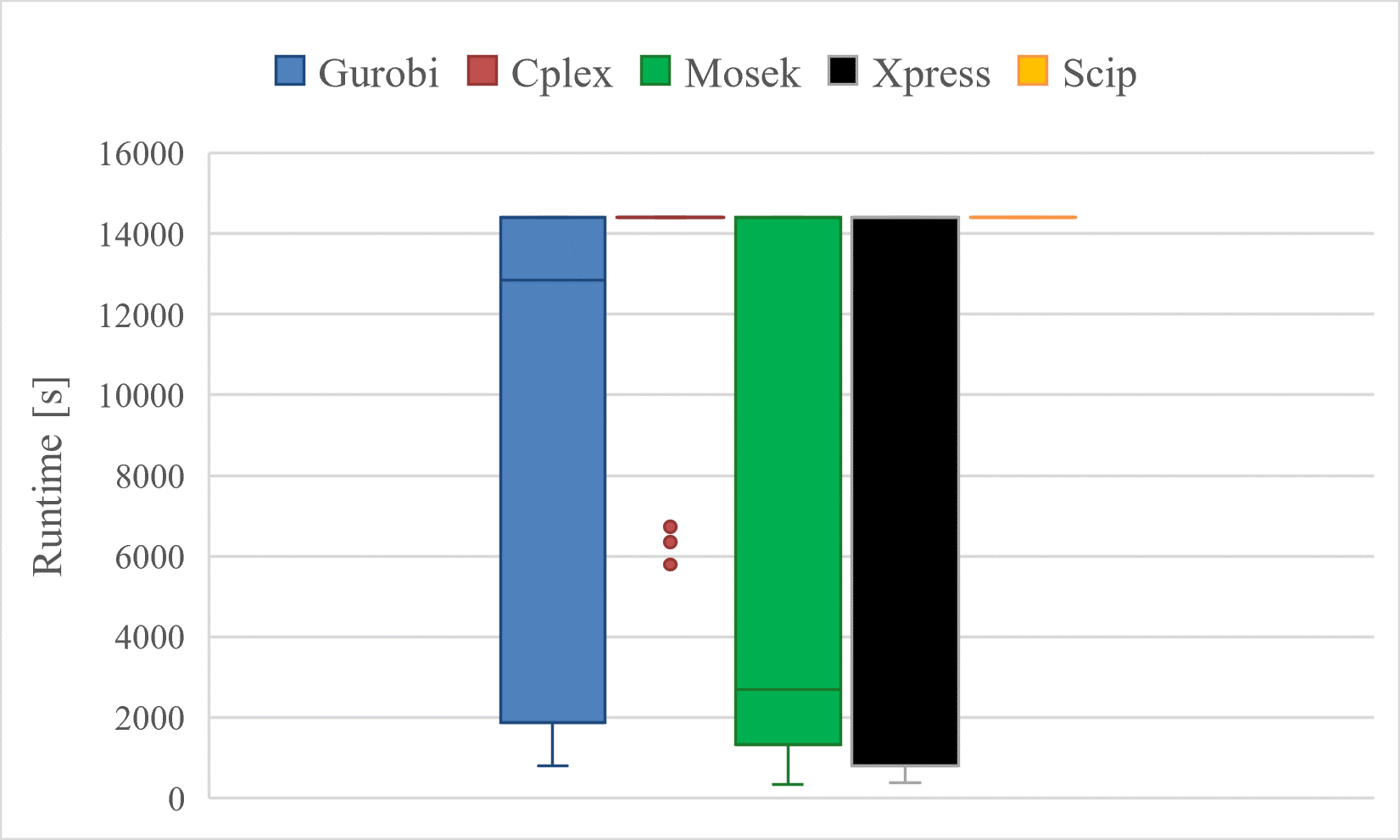}
    \caption{Runtime box-plots: small instances on the left, medium to large instances on the right.}
    \label{fig:runtime_boxplot}
\end{figure}

\begin{figure}[h]
    \centering
    \includegraphics[width=0.45\linewidth]{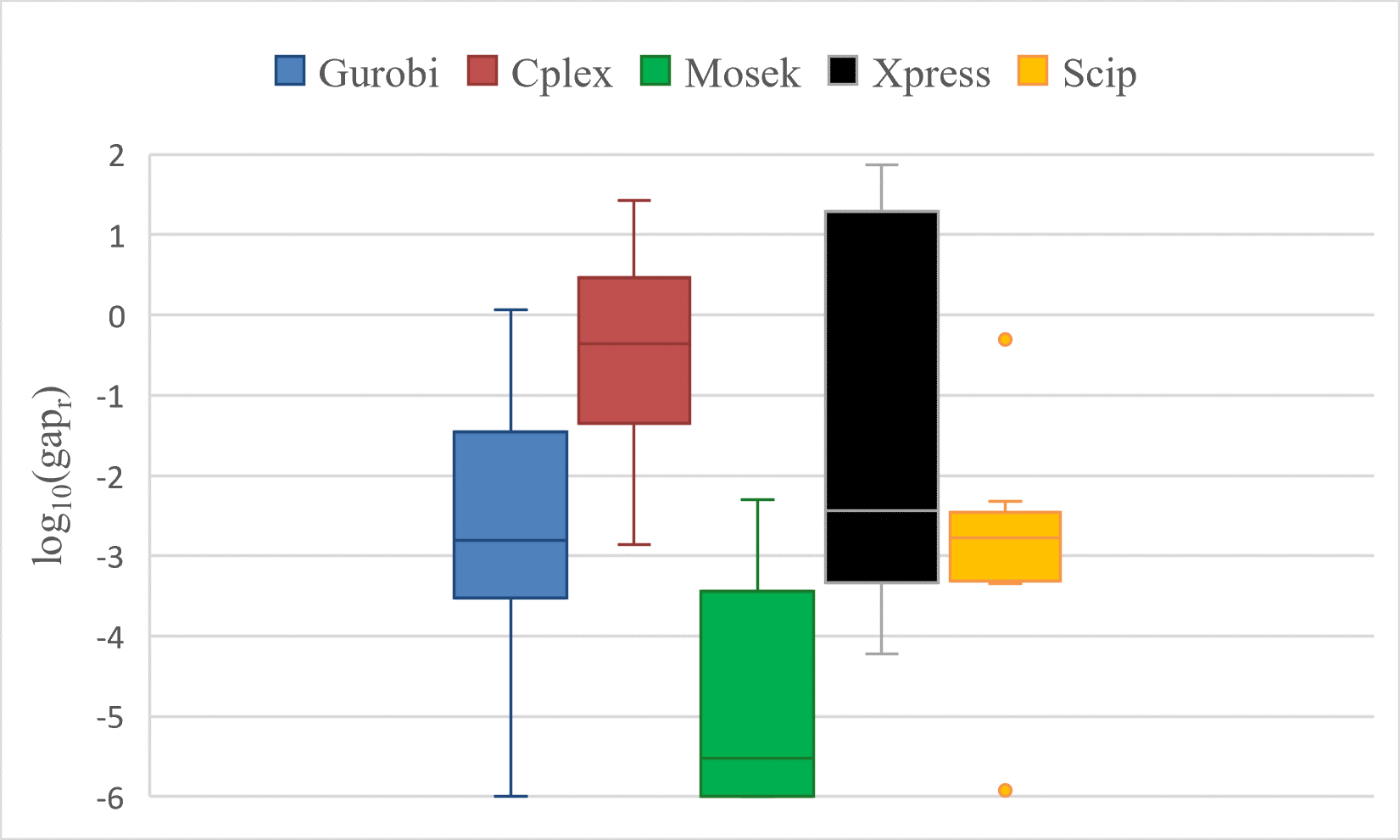}
    \includegraphics[width=0.45\linewidth]{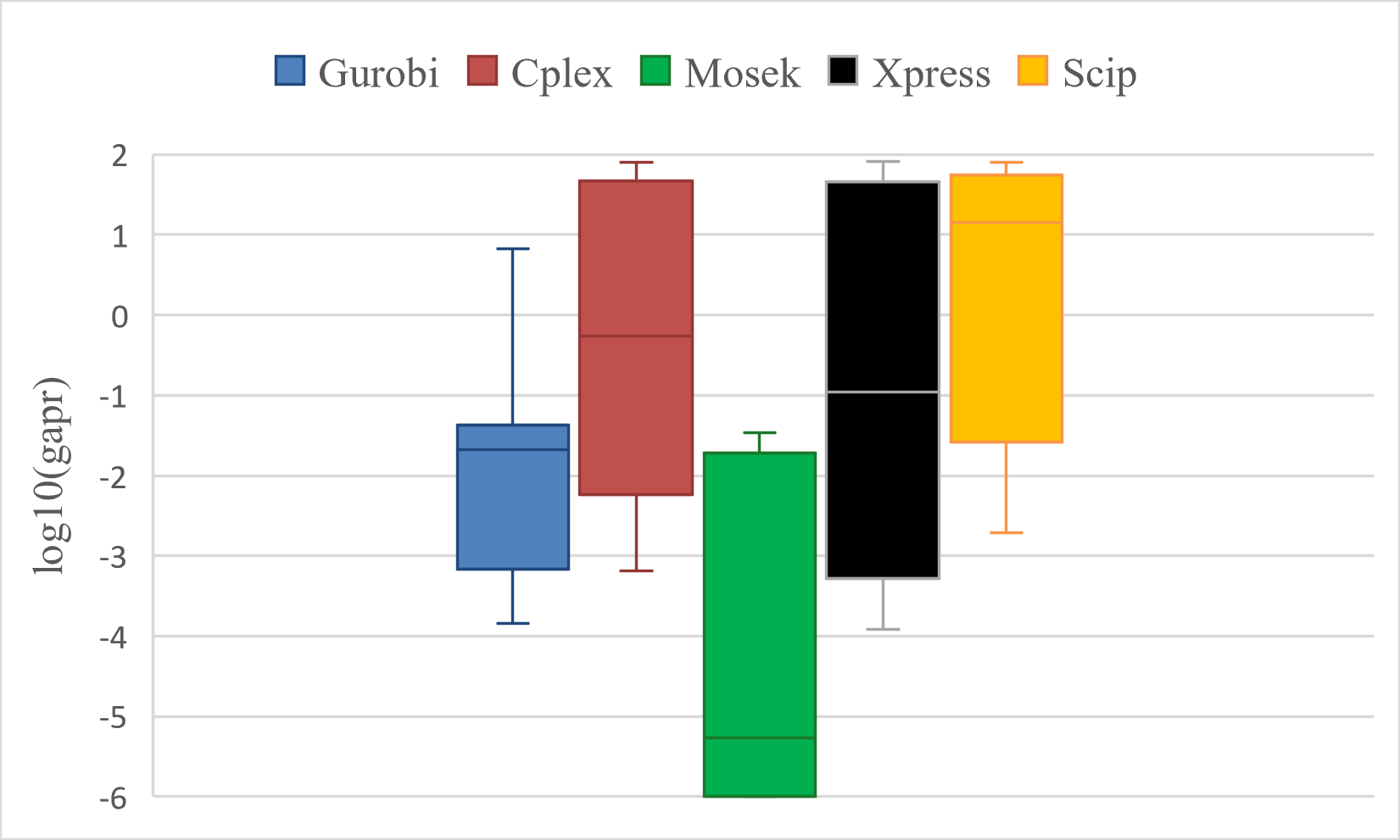}
    \caption{Logarithmic relative root gap box-plots: small instances on the left, medium to large instances on the right.}
    \label{fig:rootgap_boxplot}
\end{figure}

\begin{figure}[h]
    \centering
    \includegraphics[width=0.45\linewidth]{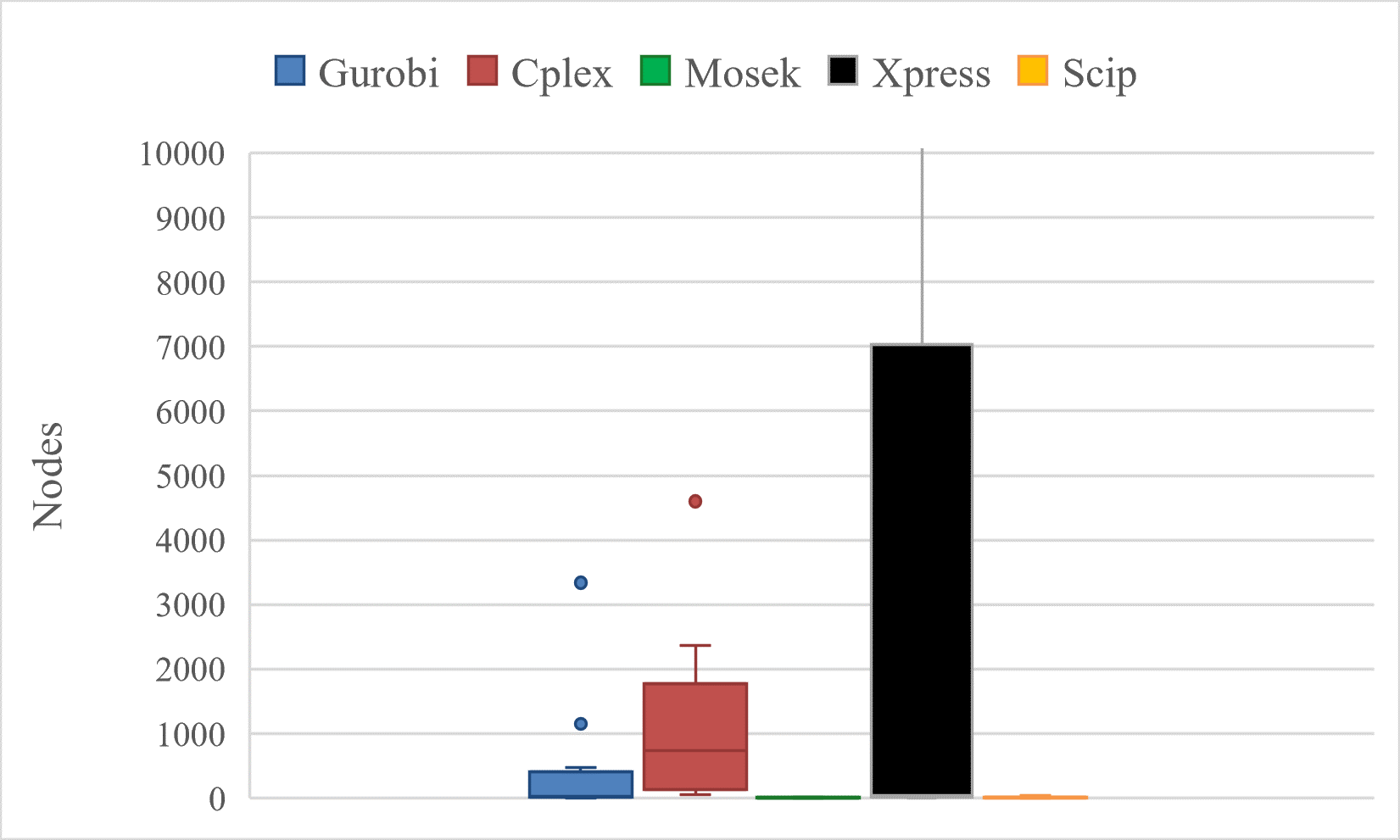}
    \includegraphics[width=0.45\linewidth]{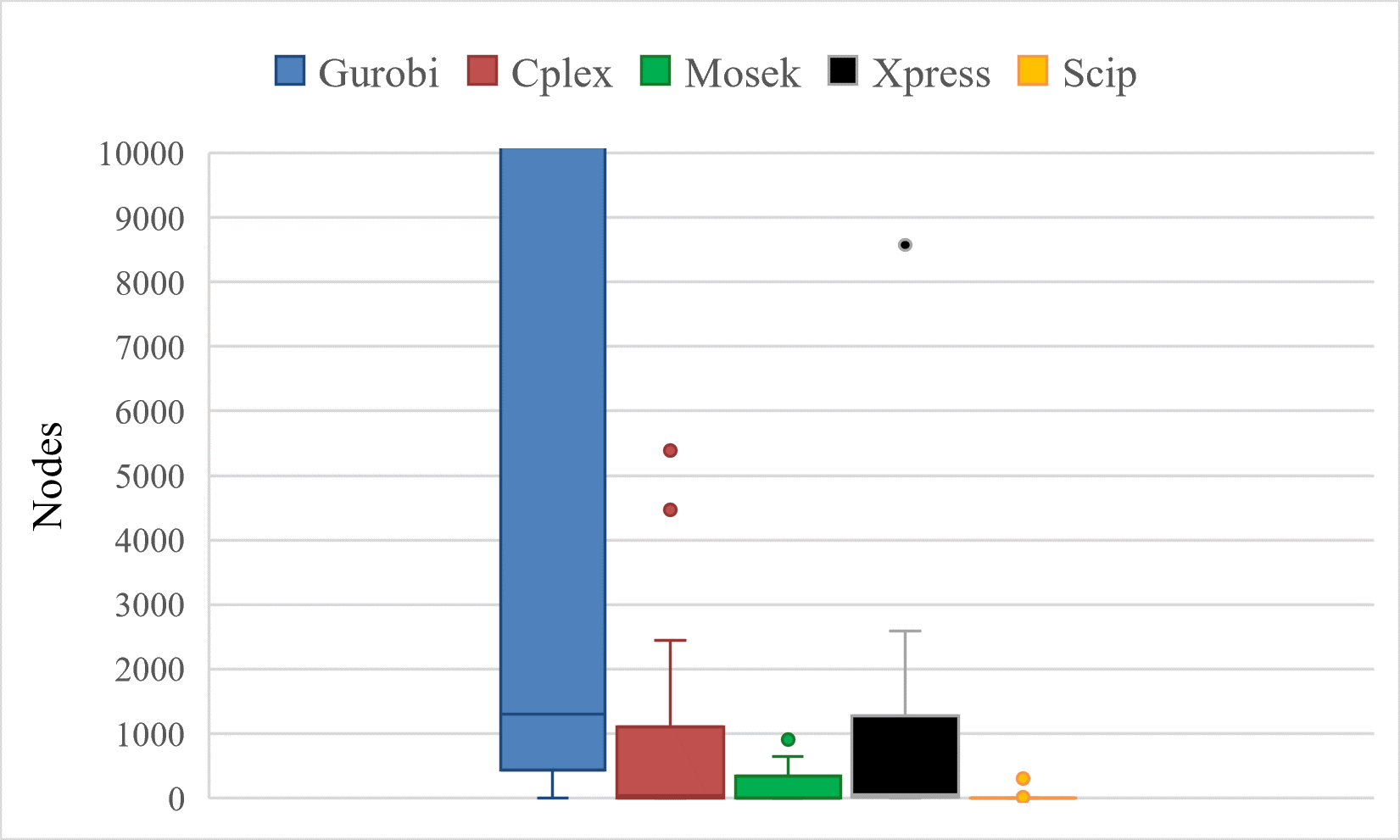}
    \caption{Number of explored nodes box-plots: small instances on the left, medium to large instances on the right.}
    \label{fig:nodes_boxplot}
\end{figure}

The second type of plot is a boxplot showing the distribution of the metric over the instances. Boxplots are constructed using the minimum value, the value of the first quartile, the median value, the value of the third quartile, and the maximum value. 
In particular, the boxes enclose the values from the first to third quartile; the horizontal line inside the box represents the median, and the two whiskers extend to the minimum and maximum values.
In Figures \ref{fig:runtime_boxplot}, \ref{fig:rootgap_boxplot} and \ref{fig:nodes_boxplot}, the boxes for the five solvers are reported. 
The median line of the box represents the central tendency of the solver's performance on the chosen metric. Additionally, the box size represents the spread of the solver's performance across the instances values, with small boxes indicating that the solver's performance is uniform.
 The plots are split between small instances (the first two blocks of instances) on the left picture, and medium to large instances (the remaining three blocks) on the right picture.

\begin{figure}[h!]
    \centering
    \includegraphics[width=0.7\linewidth]{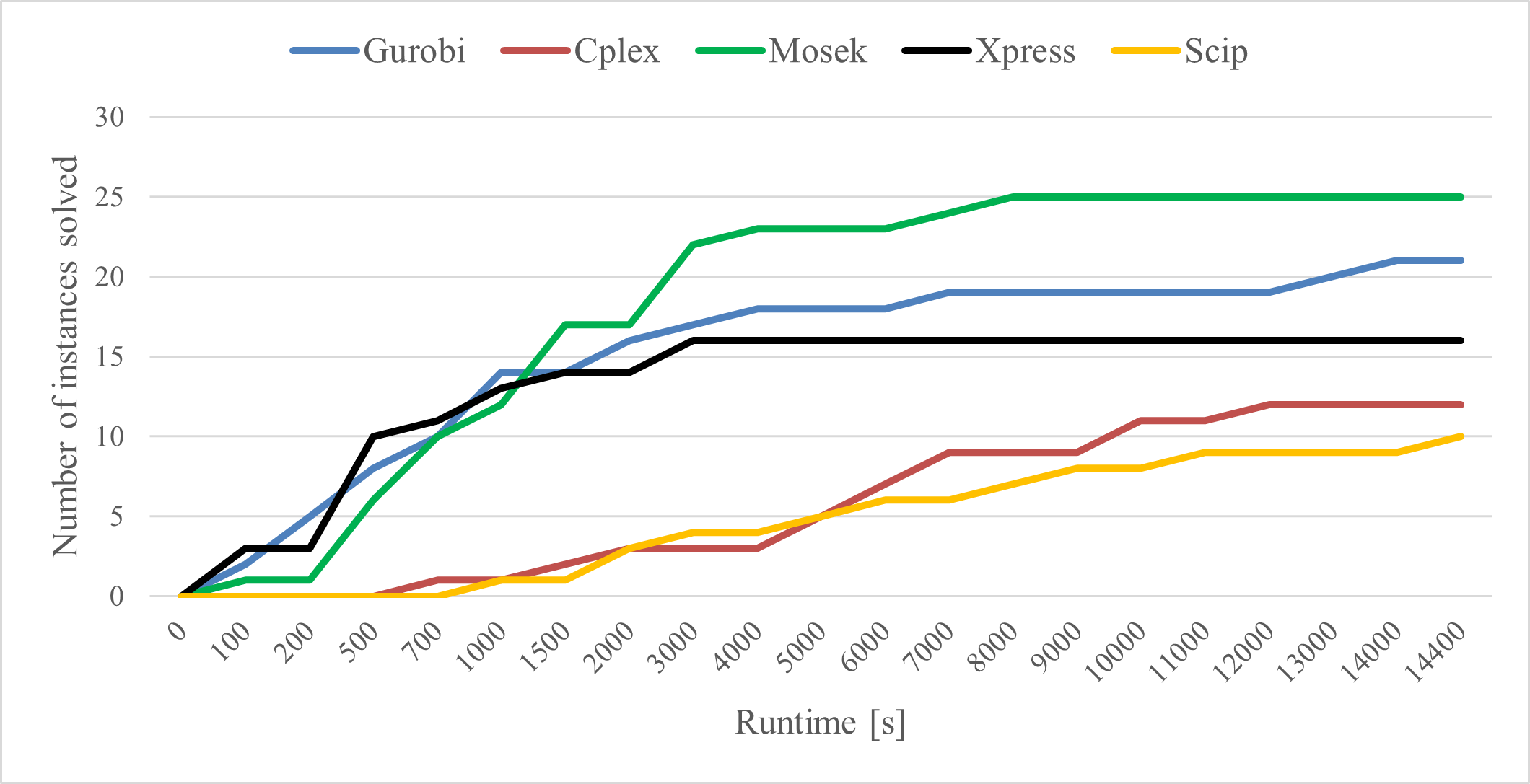}
    \caption{Solution profiles indicating the number of solved instances as a function of time.}
    \label{fig:solution_profile}
\end{figure}

We also report in Figure \ref{fig:solution_profile} the solution profiles, which represent the number of solved instances as a function of time for each solver. 
We decided not to employ performance profiles where the time is normalized with respect to best solver \cite{dolan2002benchmarking} since these are not necessarily good for comparing more than one solver simultaneously, as mentioned in \cite{gould2016note,kronqvist2019review}.

\begin{figure}[h!]
    \centering
\includegraphics[width=0.7\linewidth]{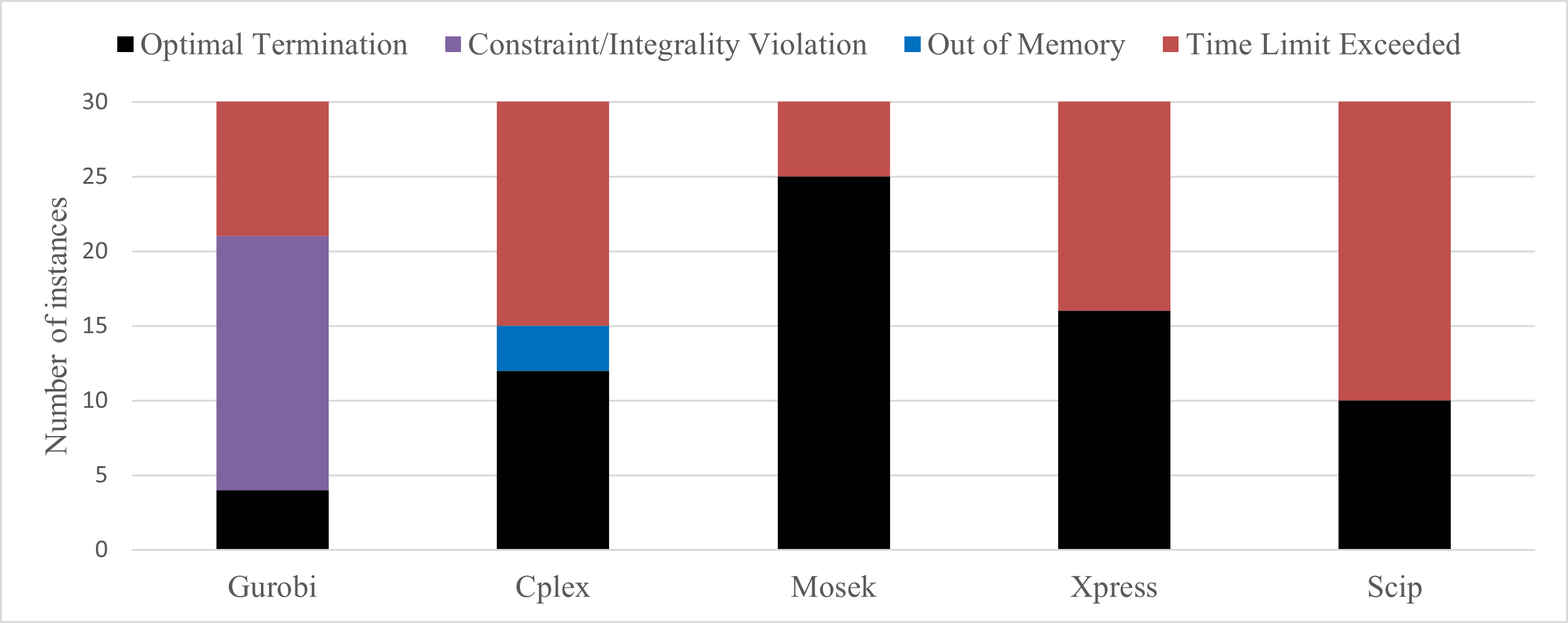}
    \caption{Solution status returned from the solvers.}
\label{fig:optimally_solved_instances}
\end{figure}

\begin{figure}[h!]
    \centering
    \includegraphics[width=0.7\linewidth]{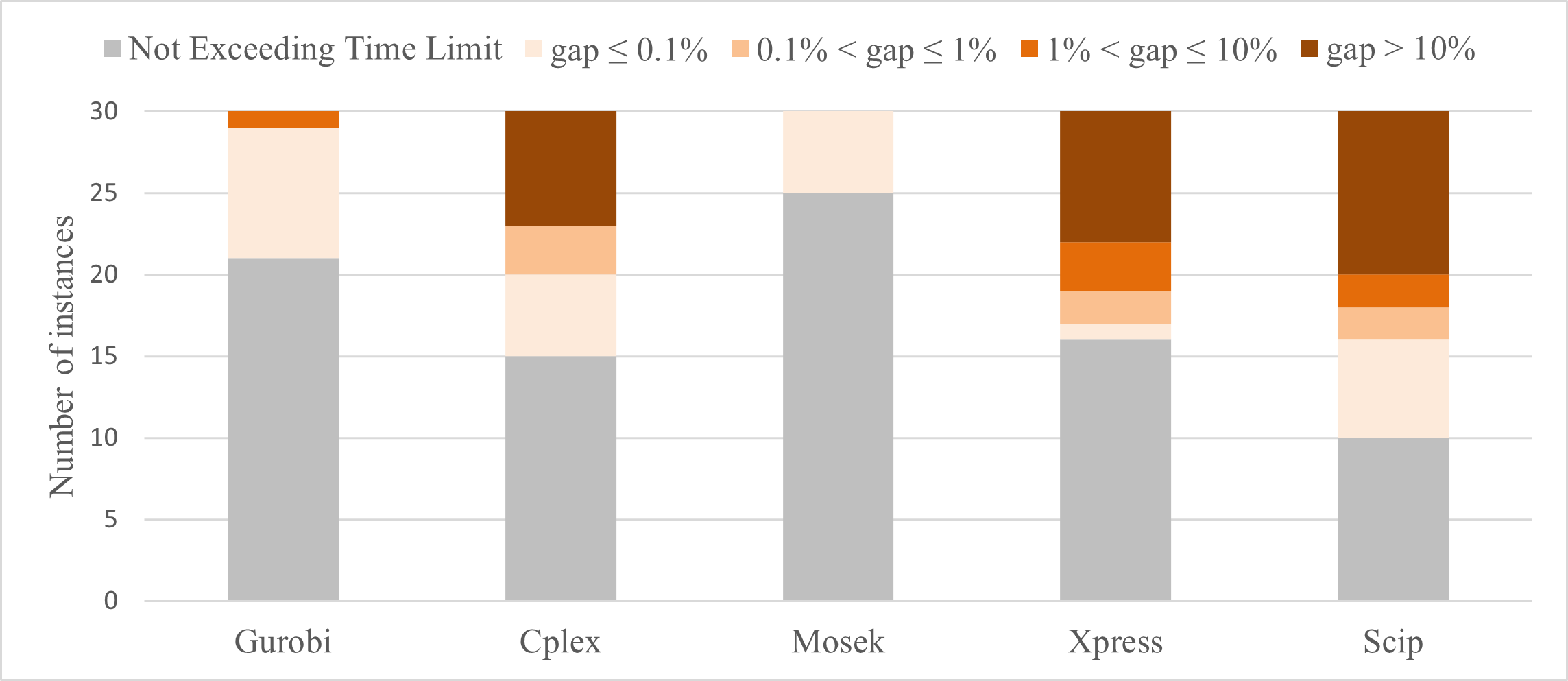}
    \caption{Relative gap at the end of the optimization when the time limit is exceeded.}
    \label{fig:optimally_solved_instances_gap}
\end{figure}

Figures  \ref{fig:optimally_solved_instances} and  \ref{fig:optimally_solved_instances_gap} display the solution status reported by the solvers and the relative gap at the end of the optimization for the instances in which the time limit was exceeded.

In the following, we discuss each metric separately.

\paragraph{Runtime}
Figures \ref{fig:runtime} and \ref{fig:runtime_boxplot} provide an illustration of the solvers' performance with regard to runtime.

Figure \ref{fig:runtime} reports on the y-axis the runtime of each solver on each instance.  The position of the symbol for each solver indicates the amount of time taken to solve the instance. 
The lower the position of the symbol, the closer the runtime is to zero and the faster the solver is in solving the instance. On the other hand, if the symbol is at the upper limit of the vertical axis, it means that the solver has reached the time limit before finding an optimal solution.


Figure \ref{fig:runtime_boxplot} provides a more detailed view of the solvers' performance through the box plots.
As the median line of the box represents the central tendency of the solver's runtime, the lower the median values, the fast performance.

From Figures \ref{fig:runtime} and \ref{fig:runtime_boxplot}, we can observe that Mosek and Gurobi are the fastest solvers, with Gurobi demonstrating faster performance on small instances, while Mosek excels on larger ones. When the solution is found within the time limit, Xpress exhibits comparable performance to Gurobi and Mosek. In some large instances, Xpress outperforms Gurobi and Mosek. On the other hand, Scip and Cplex can be considered the slowest solvers, with Scip being capable of handling small instances in medium-low times and Cplex able to solve medium-sized instances in a reasonable time.

\paragraph{Root gap}

The performance of the solvers in terms of the relative gap at the end of root node analysis is displayed in Figures \ref{fig:rootgap} and \ref{fig:rootgap_boxplot}. 
By root node analysis we denote all the operations done at the root node, such as applying heuristics and generating cuts.
Given the optimal value $opt$ (i.e. the best upper bound) and the lower bound (i.e., the relaxation value) at the end of the root node analysis $LB_r$, the gap $gap_r$ is computed as $gap_r=100\frac{opt - LB_r}{opt}$. Since relative root gap values are very small, we use a logarithmic scale (base 10) for their representation; the value $-6$ in the logarithmic scale represents relative root gaps less than or equal to $10^{-6}$. Positive logarithmic values indicate that the relative root gap is greater than 1\%.

The lower the solver symbol is displayed in Figure \ref{fig:rootgap}, the closer to zero the relative root gap reached by the solver in the instance. In Figure \ref{fig:rootgap_boxplot}, the lower the median values, the closer to zero the root gap is.
\\In most cases, Mosek provides close-to-optimal or even optimal solutions at the root node. A branch-and-bound method based on a QCP solver, which typically results in very tight bounds, is primarily responsible for this outcome. Mosek is closely followed by Gurobi, which produces tight bounds across all instances. Although Scip performs similarly to Gurobi in small instances, its performance worsens in medium to large instances. Xpress also delivers comparable outcomes to Gurobi in half of the instances; however, it reaches pretty large gaps in the other half. Instance size does not appear to impact Xpress performance. Cplex achieves relatively high root gap values, with a median across all instances of around 1\%.
It is worth noting that although Gurobi, Scip, and Cplex all start from the same initial relaxation value at the root node, Gurobi obtains values much closer to the optimum at the end of the root node analysis. This suggests that the root node operations  of Gurobi, such as applying heuristics and generating cuts, are more effective than those of the other solvers.

\paragraph{Nodes}
Figures \ref{fig:nodes} and \ref{fig:nodes_boxplot} show the performance of the solvers with regard to the number of branching nodes. We only plot values in the range $[0, 10\,000]$ to avoid crushing the graph due to scaling factors. Whenever the symbol of a solver does not appear in correspondence to an instance, this means that the solver exceeded the 10\,000 nodes.

The plots of the explored nodes should be analyzed in combination with the runtime and relative root gap plots. Mosek and Scip have similar box plots for the number of branching nodes, which are much smaller with respect to the other three solvers. However, the two solvers exhibit
entirely different behaviors in terms of runtime required and root gap, as a result of their different methodological approaches. Mosek explores a small number of nodes and also requires a short time to solve the instances, demonstrating the effectiveness of its IP methods in finding tight bounds since the beginning. This is evident from the runtime and relative root gap plots.
On the other hand, Scip visits a small number of nodes similarly. However, it cannot solve most problems to optimality within the time limit, as indicated by the runtime plots, as it spends a significant amount of time at the root node attempting to reduce the gap. Although using LP solvers results in very fast node processing, Scip takes longer to explore nodes due to its focus on applying heuristics to improve the best upper bound. This can be deduced from the log generated by Scip. Thus, we can infer that Mosek has superior performance with respect to Scip. Regarding the other solvers, Gurobi, Cplex, and Xpress generally explore several more nodes than Mosek.

\paragraph{Solution profiles}
We now discuss the solution profiles in Figure \ref{fig:solution_profile}, illustrating the number of problems solved by each solver as a function of time. It is important to note that the profiles do not represent the cumulative solution time but show how many distinct problems the solvers can solve within a certain amount of time. 

According to Figure \ref{fig:solution_profile}, we can rank the solvers based on the number of problems solved within a specific time. The top and least-performing solvers are depicted at the highest and lowest lines in the figure, respectively. A large portion of the test instances is manageable for at least two solvers (Mosek and Gurobi), while the test set is challenging for two others (Cplex and Scip). Mosek and Gurobi solved most instances within the time limit (14\,400 s). Xpress solved almost half of the instances. We also note that Xpress solved one-third of the instances in less than 500 seconds. The highest-performing solver (Mosek) solved 25 problems out of the 30 test problems, while the lowest-performing solver (Scip) only solved 10. 
None of the solvers can solve all 30 instances.
This vast disparity in the solvers' performance highlights the importance of selecting software well-suited for the specific problem type.

\paragraph{Termination status}
Lastly, we report the statistics regarding the termination of the solvers in Figure \ref{fig:optimally_solved_instances} and the statistics on the relative gap at the end of the optimization for the instances in which the time limit was exceeded in Figure \ref{fig:optimally_solved_instances_gap}. 

Mosek solved 83\% of the problems to optimality, followed by Xpress (53\%), Cplex (40\%), and Scip (33\%). It is worth noting that Gurobi declared some solutions as optimal, which actually violated the imposed tolerances on feasibility and/or integrality. Such violations are reported as a warning in the Gurobi log file. In some cases, they were significant (the maximum violation is of the order of $10^{- 1}$). The same issue was also mentioned in \cite{socp_benchmark}. 

Consequently, it is important to distinguish between the solutions that are truly optimal (13\%) and those that violate some constraints or integrality requirements (57\%).

If we consider only the instances that could not be solved within the time limit, we can differentiate the solvers based on the relative gaps reached at the end of the optimization, as shown in Figure \ref{fig:optimally_solved_instances_gap}. The relative gap at the end of the optimization is computed as $100\frac{opt - LB}{opt}$, where $opt$ is the optimal value (i.e. the best upper bound) and $LB$ is the best lower bound. As illustrated in the figure, Gurobi and Mosek achieve relative gaps of less than 0.1\% in most cases, while Cplex, Xpress, and Scip exhibit relative gaps of greater than 10\% in almost half of the instances considered.

\section{Conclusions}
\label{ch:conclusions} 

Our study aimed to evaluate the performance of Mosek, Gurobi, Xpress, Cplex, and Scip on the congested capacitated facility location problem. We analyzed the solution strategies adopted by the solvers and compared their computational efficiency on a set of 30 instances taken from the literature, in which the problem is formulated as a mixed-integer conic program.

The results revealed that Mosek and Gurobi offer competitive optimization performance, with Mosek outperforming Gurobi in large-size problems and providing more accurate solutions in terms of feasibility. Although the number of problems solved by Xpress within the time limit is less than those of Mosek and Gurobi, Xpress has demonstrated performance comparable to that of Gurobi and Mosek in half of the instances. On the other hand, Cplex and Scip emerged as the least competitive solvers on this specific testbed, with Scip's performance being affected by the default choice of spending a significant amount of time searching for a good feasible solution already at the root node. We remark that Scip is the only non-commercial solver among the ones we tested.

The results of this study highlight the importance of choosing an appropriate optimization solver for the specific problem type. Mosek's success can be attributed to its close estimate of the bounds already at the root node (it solved almost 57\% of the instances at the root), making it a suitable choice for large-size CCFLP problems. Meanwhile, for smaller problems, both Mosek and Xpress may be viable alternatives as they are more accurate and reliable than Gurobi in terms of feasibility of the declared optimal solution.

\section*{Acknowledgments}

The authors wish to thank Matteo Fischetti, Ivana Ljubi\'c, and Markus Sinnl, for providing the test instances used in this study.

\bibliographystyle{unsrt}

\end{document}